\documentclass[twocolumn,10pt]{article}

\usepackage[dvips]{graphicx}

\begin{document}

\twocolumn[

\begin{center}
{\LARGE \bf Minimum Variance Solution of Underdetermined} \\
\vspace{10pt}

{\LARGE \bf  Systems of Linear Equations} \\

\vspace{10pt}
Lorenzo Piazzo$^a$, Davide Elia$^b$, Sergio Molinari$^b$ \\

\vspace{5pt}
a) DIET, Sapienza University of Rome, Italy\\
b) IAPS, Istituto Nazionale di Astrofisica, Rome, Italy

\vspace{5pt}
E-mail: lorenzo.piazzo@uniroma1.it  \\
DIET dept. Tech. Note, June 15th, 2019 \\

\end{center}

]


{\bf Abstract.} {\em A system of linear equations is said underdetermined when there are more unknowns than equations. Such systems may have infinitely many solutions. In this case, it is important to single out solutions possessing special features. A well known example is the Minimum Norm (MN) solution, which is the solution having the least Euclidean norm. In this note, we consider another useful solution, related with the MN one, which we call the Minimum Variance (MV) solution. We discuss some of its properties and derive a simple, closed form expression.
}

\section{Introduction}

A system of linear equations is said underdetermined when there are more unknowns than equations. Such systems are found in several practical applications, an example being the estimation of $n$ parameters of a physical system from a set of $m < n$ observations. An underdetermined system can have no solution or infinitely many solutions. In the latter case, the problem of selecting a solution among all the possible ones arises. A classical way to tackle this  problem is to assign a cost to each solution and to pick the one having the minimum cost \cite{madych}. The cost function shall be selected based on the specific application and shall reflect some properties of the underlying physical system. However, there are some choices which are suitable for many applications, because they reflect general requirements, that most physical systems verify. A well known example is the Minimum Norm (MN) solution, which is the solution having the minimum Euclidean norm. A second example is the Maximum Entropy solution.

In this note, we consider another useful solution, related with the MN one, which we call the Minimum Variance (MV) solution. As we will see, while the MN solution is as close as possible to zero, the MV solution is as close as possible to its mean. As a result, the MV solution is usually smoother than the MN one, which is a desirable feature in many applications. The MV solution can be obtained numerically, by solving a convex optimization problem. Clearly, a closed form expression would be useful but, to the best of our knowledge\footnote{We scanned the Internet, chatted with colleagues and checked some textbooks \cite{meyer,lawson}. Perhaps, a closed form expression exists but it is buried under a cryptic mathematical notation. In the latter case, this note may still be useful, since we use only basic linear algebra concepts.}, it cannot be found in the literature. Such an expression is presented in this note.

{\bf Organization.} In section \ref{sec_pre} we give some definitions and basic facts. In section \ref{sec_min} we present the main result, namely a closed form for the MV solution. In section \ref{sec_der} we present the derivation.

{\bf Notation.} We use lowercase letters to denote vectors and uppercase letters to denote matrices, e.g. $a$, $B$. We denote the elements using a subscript, $a_i$, $B_{i,j}$. We use a parenthesis to specify the elements of a row vector, $a = ( a_1, ..., a_n)$. We use a superscript $T$ to denote matrix or vector transposition, $a^T$, $B^T$. We denote by $u$ the constant column vector, with all unitary elements, that is $u = ( 1, 1, ..., 1 )^T$.

\section{Preliminaries}
\label{sec_pre}

Consider a real\footnote{For the sake of simplicity, we work with real numbers. The extension to the complex case is straightforward.}, $n \times 1$ vector $x = ( x_1, ..., x_n)^T$. Its squared Euclidean norm is denoted by $|x|^2$ and is given by
\[
|x|^2 = x^T x = \sum_{i=1}^{n} x_i^2.
\]
The norm of $x$ is a measure of how far the elements of $x$ are from zero. The arithmetic mean of $x$ is denoted by $\mu_x$ and is 
\[
\mu_x = \frac{u^T x}{ u^T u } = \frac{ 1}{ n } \sum_{i=1}^{n} x_i .
\]
The vector $x$ can be decomposed as follows
\[
x = \hat{x} + \bar{x} ,
\]
where
\[
\bar{x} = u \mu_x 
\]
is termed the constant part of vector $x$ and is a $n \times 1$ constant vector, with elements equal to the mean, $\bar{x} = ( \mu_x, ..., \mu_x)^T$, and 
\[
\hat{x} = x - \bar{x} 
\]
is termed the variable or zero mean part of vector $x$ and is a zero mean vector. 

With a slight abuse of terminology, the squared norm of the variable part, namely
\[
|\hat{x}|^2 = \sum_{i=1}^{n} (x_i - \mu_x)^2,
\]
will be called the arithmetic variance\footnote{Note that if the elements of $x$ are drawn from a random distribution, $|\hat{x}|^2 / n$ is an estimator of the distribution's variance.} of vector $x$. The variance of $x$ is a measure of how far the elements of $x$ are from the mean of $x$. Note that the variable and constant parts are orthogonal vectors\footnote{Two vectors $v$ and $w$ are said orthogonal if their product is zero, $v^T w = 0$.}. Therefore
\[
|x|^2 = |\hat{x}|^2 + |\bar{x}|^2 .
\]
%

\section{Minimum Variance solution}
\label{sec_min}

Consider a linear system of $m$ equations and $n$ unknowns, written in matrix form as 
\begin{equation}
A x = b
\label{sys}
\end{equation}
where $A$ is an $m \times n$ coefficient matrix, $x$ is an $n \times 1$ unknown vector and $b$ is an $m \times 1$ known vector. We assume that $m < n$, so that the system is underdetermined. We further assume that $A$ is full rank, so that $AA^T$ is non singular and the system has infinitely many solutions \cite{meyer}. In order to select a solution, a classical approach is that of assigning a scalar, real cost function, $f(x)$, and seeking the solution that minimizes $f(x)$. In other words, we have to solve the following constrained optimization problem 
\[
\begin{tabular}{l} 
minimise  $f(x)$ \\
such that $Ax = b$ .
\end{tabular}
\]
The cost function is selected based on the application. Depending on the function, the optimization problem may be hard or simple to solve numerically. Moreover, for some simple functions, we can derive a closed form solution.

An important example is when the cost function is the squared Euclidean norm, i.e. $f(x) = |x|^2$. In this case we obtain the MN solution, denoted by $x_N$, for which a closed form exists, e.g. \cite{madych,meyer,lawson}, and is
\begin{equation}
 x_N = A^T ( A A^T)^{-1} b.
\label{mns}
\end{equation}
Note that the MN solution is unique, i.e. if $z$ is a solution and $z \neq x_N$, then $|x_N|^2 < |z|^2$.
 
As another option, we can set the cost equal to the variance, i.e. $f(x) = | \hat{x} |^2$. In this case, we obtain the MV solution. As we mentioned, while the MN approach pushes the solution towards zero, the MV approach pushes the solution towards its mean, which is a desirable feature for many applications. 

The MV solution can be obtained numerically, by solving a convex optimization problem. Even better, we can derive a closed form expression for it. Specifically, in the next section we show that when $Au \neq 0$, the MV solution, denoted by $x_V$, is unique and is 
\begin{equation}
 x_V = A^T ( A A^T)^{-1} (b   -  A u \alpha ) + u \alpha,
\label{mvs}
\end{equation}
where the scalar $\alpha$ is
\begin{equation}
\alpha = \frac{  u^T A^T ( A A^T)^{-1} b       }{    u^T A^T ( A A^T)^{-1} A u         }   .
\label{alpha}
\end{equation}
Moreover, when $Au = 0$, the MV solution is not unique and the MN solution is a MV solution too.

\section{Derivation}
\label{sec_der}

Unless otherwise stated, we assume that $A$ is full rank, so that $AA^T$ is a non singular, symmetric matrix \cite{meyer}. We also assume that $Au \neq 0$. To proceed, we introduce the vector $h = Au$ and consider the following system
\begin{equation}
A x = b - h \alpha.
\label{asys}
\end{equation}
The latter system will be said a system associated to the original system (\ref{sys}), with parameter $\alpha$. Note that, if $s$ is a solution of the original system, then $v = s - u \alpha$ is a solution of the associated system. Indeed
\[
Av = A (s - u \alpha ) = A s - A u \alpha = b - h \alpha.
\]
Moreover, note that if a system is associated to a second system with parameter $\alpha$, then the second system is associated to the first one with parameter $-\alpha$. 

By varying $\alpha$ over the real axis, an infinite number of associated systems is produced. Using (\ref{mns}), for each of these systems we can compute a unique MN solution, denoted by $x_N(\alpha)$, given by
\begin{equation}
 x_N(\alpha) = A^T ( A A^T)^{-1} (b - h \alpha).
\label{mnsa}
\end{equation}
The following Lemma identifies a special, important associated system.

\vspace{5pt}
\noindent
{\bf Lemma 1. Base system and solution.} Consider the set of systems (\ref{asys}) for $-\infty < \alpha < \infty$ and the corresponding MN solutions of (\ref{mnsa}). When  $\alpha = \alpha^*$, with
\begin{equation}
\alpha^* = \frac{  u^T A^T ( A A^T)^{-1} b       }{    u^T A^T ( A A^T)^{-1} A u         }   ,
\label{astar}
\end{equation}
the resulting system will be called the base system. The corresponding MN solution, namely $x_N( \alpha^* )$, will be denoted by $x_B$ and called the base solution. It has the following properties:
\begin{equation}
| x_B |^2  < | x_N(\alpha) |^2     \hspace{.3cm} \mbox{ for }  \alpha \neq \alpha^*
\label{p1}
\end{equation}
and 
\begin{equation}
x_B  \hspace{.3cm} \mbox{ is a zero mean vector }.
\label{p2}
\end{equation}

\noindent
{\bf Proof.} By introducing the vector $p = A^T ( A A^T)^{-1} h = A^T ( A A^T)^{-1} A u$ and by denoting by $x_N$ the MN solution of the original system, given by (\ref{mns}), we can write the MN solutions of (\ref{mnsa}) in a compact form, namely as
\[
x_N( \alpha ) = x_N - p \; \alpha.
\]
The squared norm is
\[
| x_N( \alpha ) |^2 = | x_N |^2 + | p |^2 \alpha^2 - 2 p^T x_N \; \alpha .
\]
As a function of $\alpha$, the latter expression is quadratic and has a single minimum, which can be found by setting the derivative to zero. Doing so, we get
\[
2 | p |^2 \alpha - 2 p^T x_N = 0 ,
\]
and solving for $\alpha$ we have
\[
\alpha^* =  \frac{        p^T x_N       }{     | p |^2     } .
\]
By writing $\alpha^*$ explicitly we have
\[
\alpha^* =  \frac{     h^T ( A A^T)^{-1}  A A^T ( A A^T)^{-1} b      }{      h^T ( A A^T)^{-1}  A A^T ( A A^T)^{-1} h     }  = 
\]
\[
 =    \frac{     u^T A^T ( A A^T)^{-1} b      }{      u^T A^T ( A A^T)^{-1} A u     }  
\]
which proves (\ref{p1}). 

The base solution can be written explicitly as
\begin{equation}
 x_B = A^T ( A A^T)^{-1} (b - h \alpha^*) = x_N - p \alpha^*.
\label{xb}
\end{equation}
In order to show that $x_B$ is zero mean, we need to show that $u^T x_B = 0$. Indeed
\[
u^T x_B = u^T A^T ( A A^T)^{-1} b    -    u^T A^T ( A A^T)^{-1} A u   \alpha^*  = 
\]
\[
 = u^T A^T ( A A^T)^{-1} b    -   u^T A^T ( A A^T)^{-1} b = 0.
\]
$\bullet$

\vspace{5pt}

As we have seen, $x_B$ is the MN solution of the base system. As we show in the next Lemma, it is the MV solution too.

\vspace{5pt}
\noindent
{\bf Lemma 2.} The base solution $x_B$ is the unique MV solution of the base system.

\noindent
{\bf Proof.}  Suppose there is a vector $z = \hat{z} + \bar{z}$, with $z \neq x_B$, which is the MV solution of the base system. This implies that 
\[
| \hat{z} |^2  \leq  | \hat{x}_B |^2 
\]
and, since $x_B$ is zero mean, that 
\begin{equation}
| \hat{z} |^2  \leq  | x_B |^2 . 
\label{c2}
\end{equation}
We show that such a vector does not exist. We separate the proof in two cases, depending on whether $z$ is zero mean or not.

If $z$ is zero mean we have $z = \hat{z}$. Replacing in (\ref{c2}) we get
\[
| z |^2  \leq  | x_B |^2 .
\]
However, $|z|^2$ cannot be less than $|x_B|^2$, because $x_B$ is the MN solution. Neither $|z|^2$ can be equal to $|x_B|^2$ and $ z \neq x_B $, because the MN solution is unique.

If $z$ has mean $\mu \neq 0$, it can be written as $z = \hat{z} - u \mu$. In this case, $\hat{z}$ is the solution of an associated system with parameter $\beta = \alpha^* + \mu \neq \alpha^*$. Moreover, from (\ref{c2}), its squared norm is less than or equal to $|x_B|^2$. But the latter facts contradict (\ref{p1}). $\bullet$

\vspace{5pt}

We are now ready to prove our main result. First, we note that $s = x_B + u \alpha^*$ is a solution of the original system and is in fact the solution given by (\ref{mvs}), as is easy to check. In the next Lemma, we show that $s$ is indeed the MV solution of the original system.

\vspace{5pt}
\noindent
{\bf Lemma 3.} The solution $s = x_B + u \alpha^*$ is the unique MV solution of the system (\ref{sys}).

\noindent
{\bf Proof.} Assume that the original system has a different MV solution, denoted by $z$. Note that the variance of $s$ is $|x_B|^2$, therefore, in order to be the MV solution, the variance of $z$ has to be less than or equal to $|x_B|^2$. However, if $z$ had a variance lower than $|x_B|^2$, by adding a proper constant to $z$, we could construct a solution for the base system having variance lower than $|x_B|^2$, which contradicts Lemma 2. Then the variance of $z$ must be  equal to that of $s$, i.e. $|\hat{z}|^2 = |x_B|^2$.

Now, denote by $\mu$ the mean of $z$ so that $z = \hat{z} + u \mu$ and consider the following associated system 
\[
Ax = b - h \mu .
\]
Note that $\hat{z}$ is a solution of the latter system. Moreover, its squared norm is $|x_B|^2$. Then, from (\ref{p1}), we have that $\mu = \alpha^*$, showing that the latter system is in fact the base system and that $\hat{z}$ is a MN solution of the base system, because $|\hat{z}|^2 = |x_B|^2$. Since the MN solution is unique, we conclude that $\hat{z} = x_B$. $\bullet$

\vspace{5pt}

To conclude, we consider the case when $Au = 0$. In this case, given any solution $s$, we can add or subtract a constant vector to $s$ and obtain another solution. Since the constant part does not affect the variance, the latter fact implies that the MV solution is not unique. Moreover, the MN solution is zero mean, otherwise we could reduce its norm by subtracting the mean. Furthermore, the minimum attainable variance is equal the minimum attainable squared norm. As a result, the MN solution is a MV solution too.


\end{document}